\documentclass{amsart}
\usepackage{verbatim}
\usepackage{graphicx}
 \newtheorem{thm}{Theorem}[section] %was [subsection]
 \newtheorem{cor}[thm]{Corollary}
 \newtheorem{lem}[thm]{Lemma}
 
 \theoremstyle{definition}
 \newtheorem{defn}[thm]{Definition}
 \theoremstyle{remark}
 
 \numberwithin{equation}{section} %was {subsection}
\begin{document}

\title[Probabilistic Inductive Classes of Graphs]
 {Probabilistic Inductive Classes of Graphs}

\author{Nata\v{s}a Kej\v{z}ar$^{1}$}
\author{Zoran Nikoloski$^{2}$}
\author{Vladimir Batagelj$^{3}$}

\address{$^{1}$Faculty of Social Sciences, University of Ljubljana, Slovenia}
\address{$^{2}$Department of Applied Mathematics, Faculty of Mathematics and Physics, Charles University, Prague, Czech Republic}

\address{$^{3}$Faculty of Mathematics and Physics, University of Ljubljana, Slovenia}

\email{natasa.kejzar@fdv.uni-lj.si, nikoloski@kam.mff.cuni.cz}

\thanks{This work was completed with the support of the Department of Applied Mathematics, Faculty of Mathematics and Physics, Charles University, Prague grant MSM0021620838; the DELIS-IST 001907 project of the European Comission; by the Ministry of Higher Education, Science and Technology of Slovenia, Project J1-6062; and partially by the Slovenian Research Agency, Project P5-0168.}

%\subjclass{Primary 47A15; Secondary 46A32, 47D20}


\date{December 27, 2006 (submitted).}

\dedicatory{}

%\commby{Daniel J. Rudolph}

%%% ----------------------------------------------------------------------

\begin{abstract}
Models of complex networks are generally defined as graph stochastic processes in which edges and vertices are added or deleted over time to simulate the evolution of networks. Here, we define a unifying framework --- \textit{probabilistic inductive classes of graphs} --- for formalizing and studying evolution of complex networks. Our definition of probabilistic inductive class of graphs (PICG) extends the standard notion of inductive class of graphs (ICG) by imposing a probability space. A PICG is given by: (1) class $\mathcal{B}$ of initial graphs, the \textit{basis} of PICG, (2) class $\mathcal{R}$ of generating rules, each with distinguished \textit{left element} to which the rule is applied to obtain the \textit{right element}, (3) probability distribution specifying how the initial graph is chosen from class $\mathcal{B}$, (4) probability distribution specifying how the rules from class $\mathcal{R}$ are applied, and, finally, (5) probability distribution specifying how the left elements for every rule in class $\mathcal{R}$ are chosen. We point out that many of the existing models of growing networks can be cast as PICGs. We present how the well known model of growing networks --- the preferential attachment model --- can be studied as PICG.
As an illustration we present results regarding the size, order, and degree sequence
for PICG models of connected and 2-connected graphs.
\end{abstract}

%%% ----------------------------------------------------------------------
\maketitle
%%% ----------------------------------------------------------------------

\section{Introduction}

Recent surge of empirical results pertaining to the characteristics of complex networks, from technological networks (the Internet, the power grid) to networks of social contacts and protein interactions, has prompted theoretical research in models of complex networks. Existing theoretical results have identified two mechanisms that result in synthetic networks with characteristics closely matching those of the real-world counterparts. Models of complex networks are usually defined as graph stochastic processes in which edges and vertices are added or deleted over time to simulate the evolution of networks, based on (1) preferential attachment mechanisms \cite{Barabasi99}, (2) copying mechanism \cite{Kumar00}, or a combination of the two \cite{Leskovec05}. The result of a graph stochastic process is a class of random graphs, which can also be obtained via \textit{probabilistic inductive classes of graphs} (PICG).

\section{Definition of PICG}

Eberhard \cite{Eberhard91} was the first to inductively define a class of graphs. A comprehensive survey of existing inductive definitions for classes of graphs is presented in \cite{Batagelj85}.

\begin{defn}
An \textit{inductive class of graphs} (\textit{ICG}) $\mathcal{I}=(\mathcal{B};\mathcal{R})$ is defined by \cite{Curry63}:
 \begin{enumerate}
 \item class $\mathcal{B}$ of initial graphs, the \textit{basis} of ICG,
 \item class $\mathcal{R}$ of generating rules, each with distinguished \textit{left element} (part of a graph) to which the rule is applied to obtain the \textit{right element}.
 \end{enumerate}
\end{defn}

\begin{figure*}
   \includegraphics[width=\textwidth]{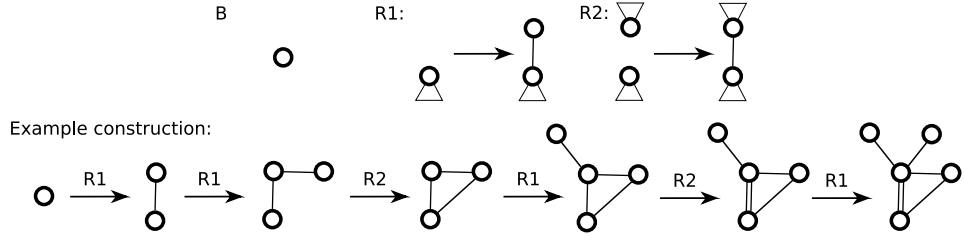}
\caption{\label{f:example} An example of an ICG class $\mathcal{I}=(B;R1,R2)$ and a possible construction of a graph in that class.}
\end{figure*}

The rules describe a substitution of a part of a graph (left side) by another part of a graph (right side) and result into a graph again.
A class $\mathcal{I}$ consists exactly of the graphs that can be constructed from the  graphs in the basis $\mathcal{B}$ by applying a finite number of generating rules $\mathcal{R}$.
Figure \ref{f:example} shows an example of an ICG $\mathcal{I}=(\mathcal{B};\mathcal{R})$, $\mathcal{R}=\{ R1,R2 \}$ and a possible construction for one of its graphs.

%The small triangles in the parts of graphs used in the pictures describing the rules denote 0 or more edges.
The small triangles attached to the nodes in the pictorial depiction of the rules presented in Fig. \ref{f:example} denote 0 or more edges.

Here we extend the classical notion of inductive class of graphs by introducing a probabilistic space in the following manner:

\begin{defn}
A \textit{probabilistic inductive class of graphs} (\textit{PICG}), $\mathcal{I}$, is given by:
\begin{enumerate}
\item class $\mathcal{B}$ of initial graphs, the \textit{basis} of PICG,
\item class $\mathcal{R}$ of generating rules, each with distinguished \textit{left element} to which the rule is applied to obtain the \textit{right element},
\item probability distribution specifying how the initial graph is chosen from class $\mathcal{B}$,
\item probability distribution specifying how the rules from class $\mathcal{R}$ are applied, and, finally,
\item probability distribution specifying how the left elements for every rule in class $\mathcal{R}$ are chosen.
\end{enumerate}
\end{defn}

The probabilistic inductive class $\mathcal{I}$ is closed under the class of generating rules $\mathcal{R}$, \textit{i.e.} a random graph $G$ is in the PICG $\mathcal{I}$ if and only if it can be obtained by applying the generating rules in the given class $\mathcal{R}$ finite number of times, starting with an initial graph from the basis $\mathcal{B}$. All steps in this construction have to have positive probabilities.

It is often desired that the generating rules are simple enough not only to facilitate rigorous analysis but also to capture the mechanism for generating a class of graphs that match the properties of real-world networks. Two such properties are:

\begin{enumerate}
\item \textit{locality}: the rule is local if and only if its left element is connected.
\item \textit{expansion}: the rule is expanding if a selected property of graphs is increased (or remains the same) by the application of the rule (for example number of vertices, number of edges, girth --- the size of the shortest cycle, etc.)
\end{enumerate}

In general it can happen that the selected rule can not be applied in a current graph. In this paper we shall limit our discussion to the \textit{simple} PICG definitions in which the left element of the selected rule can always be found.

In order to prove a certain result for ICGs (as well as PICGs), the principle of \textit{inductive generalization} is often used:

\textbf{Principle: }\textit{Inductive generalization}

Given an ICG $\mathcal{I}$ and a property of graphs $Q$. If it holds that:

\begin{enumerate}
\item each graph in $\mathcal{B}$ has the property $Q$, and
\item each generating rule in $\mathcal{R}$ \textit{preserves} the property $Q$ --- if the graph on which the rule is applied, has the property $Q$, then also the resulting graph has the property $Q$,
\end{enumerate}

then all graphs from $\mathcal{I}$ have the property $Q$.

\textit{An example}: If all the elements (initial graphs $\mathcal{B}_{PICG}$ and the rules $\mathcal{R}_{PICG}$) in the definition of PICG have positive probabilities, then $\mathcal{I}_{PICG} = \mathcal{I}_{ICG}$, where $\mathcal{I}_{ICG} = \left (\mathcal{B}_{PICG};\mathcal{R}_{PICG}\right )$, because every graph from ICG has a positive probability to be created.

Notation: Given a graph $G = \left ( V , E \right )$, let us denote the number of vertices $\left | V \left ( G \right ) \right | = n$, and the number of edges $\left | E \left ( G \right ) \right | = m$; and random variables for the number of vertices with ${\bf n}$ and ${\bf m}$ for the number of edges.  Let $d \left ( u \right )$ denote the degree of vertex $u \in V \left ( G \right )$ --- number of edges with $u$ as an endpoint.

% ------------------------------------------------------------------------

\section{\label{s:PA model}An illustration: Preferential attachment model as PICG}

The model was presented by Barab\'asi and Albert in 1999 \cite{Barabasi99}.
Its idea is to capture the evolution of the network that can possibly explain the emergence of a power-law degree distribution.
The model is based on the following concepts: (1) a network grows in time (the number of vertices and number of edges grow) and that (2) a newcomming vertex is more likely to build an edge to a more ``popular" vertex in the network (a vertex with relatively large number of edges). Barab\'asi and Albert called this principle a preferential attachment (abbreviated as \textit{pa}).

The algorithm of a model consists of two processes in a time step: an addition of a new vertex and a creation of $m_{pa}$ edges from this vertex to already existing vertices, which are chosen proportional to their current degree.
To make the formulation of the model precise, consider to start the process with two vertices linked by $m_{pa}$ paralel edges. When adding a new vertex, the edges will be added one at a time (so for the second and subsequent edges the probabilities will be calculated with updated vertex degrees). This formulation has a nice property that only the $m_{pa}=1$ model needs to be analyzed \cite{Bollobas01}. To get the models with $m_{pa}>1$, vertices, added at time steps $jm_{pa}, jm_{pa}-1, \ldots, (j-1)m_{pa}+1$, should be collapsed into a single vertex $j$ (added at time step $j$).

\begin{figure}[!thbp]
   \includegraphics[width=80mm]{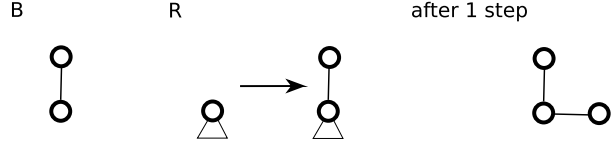}
\caption{\label{f:pa model}The basis, the only rule for a preferential attachment model and the graph after the first step (after one application of the rule).}
\end{figure}

We can describe this model as a probabilistic ICG very easily:
The basis consists of an edge linking 2 vertices. The rule is just one --- to a selected vertex an edge that ends with a new vertex is added (see Fig. \ref{f:pa model}). We only have to specify, how the left element (the vertex) is chosen. The probability of selecting a vertex, where to apply the rule, is proportional to its current degree. The probability of a vertex $i$ to be chosen is exactly $s_i =d(i)/\sum_j d(j)$.

% ------------------------------------------------------------------------

\section{Order and size of PICG}

For the simple PICG it is easy to determine the expected order (number of vertices) and expected size (number of edges) of the graph after $t$ steps. We divide the rules in three different classes: expanding, shrinking and stable class. The change in the number of elements (vertices or edges) after the application of the rule in these classes is either growing, decreasing or staying the same. In order to get non-trivial results for large $t$, the rules' expected change in the number of vertices or edges has to be positive. Otherwise the graph dies out when $t \rightarrow \infty$.

Given a graph $G \in \mathcal{I}$, let $p_t \left ( {\bf n} = n \right )$ denote the probability that at step $t$ the order of $G$ is $n$. Consider $R_i \in \mathcal{R}$ a fixed chosen rule in $\mathcal{I}$, $\Delta n_i$ the number of vertices, that rule $R_i$ adds to the graph and $r_i$ the probability of selecting rule $R_i$. The probability $p_t \left ( {\bf n} = n \right )$ can be expressed recursively as

\[
p_t \left ( {\bf n} = n \right ) = \sum_{R_i \in \mathcal{R}}{r_i p_{t-1}\left ( {\bf n} = n - \Delta n_i\right )},
\]
with initial values for all basic graphs $B_i \in \mathcal{B}$: $p_0 \left ( {\bf n} = k \right ) = \sum_{B_i;\left | B_i \right |=k} q_i$, where $q_i$ is the probability, that $B_i$ is chosen for initial graph. For the clarity of the further equations the following abreviation will be used $p_t(n) \equiv p_t ( {\bf n} = n )$.

The expected number of vertices is therefore
\[
E_t \left [ {\bf n} \right ] =  \sum_{i=0}^{\infty}{i p_t \left ( i \right )}
\]
The expected change in the number of vertices is then easily computed as
\begin{equation}
\lim_{t \rightarrow \infty}\frac{E_t \left [ {\bf n} \right ]}{t} = \sum_{R_i \in \mathcal{R}}{r_i \Delta n_i}
\label{eq:expected_change}
\end{equation}

\begin{proof}
We can write
\begin{eqnarray*}
E_t \left [ {\bf n} \right ] &=& \sum_{i=0}^{\infty} i \sum_{R_j \in \mathcal{R}} r_j p_{t-1}\left ( i - \Delta n_j\right ) \\
&=& \sum_{R_j \in \mathcal{R}} r_j \sum_{i=0}^{\infty} i p_{t-1}\left ( i - \Delta n_j\right )
\end{eqnarray*}
We split the inner sum into two sumands
\begin{eqnarray*}
\sum_{i=0}^{\infty}{i p_{t-1}\left ( i - \Delta n_j\right )}
&=& \sum_{i=0}^{\infty}{(i-\Delta n_j) p_{t-1}\left ( i - \Delta n_j\right )} + \sum_{i=0}^{\infty}{\Delta n_j p_{t-1}\left ( i - \Delta n_j\right )} \\
&=& \sum_{k=0}^{\infty}{k p_{t-1}\left ( k \right )} + \Delta n_j \\
&=& E_{t-1} \left [ {\bf n} \right ] + \Delta n_j.
\end{eqnarray*}
Now, the whole equation can be written as a recursive relation
\begin{eqnarray*}
E_t \left [ {\bf n} \right ] &=&  \sum_{R_j \in \mathcal{R}} r_j \left ( E_{t-1} \left [ {\bf n} \right ] + \Delta n_j \right ) \\
& = & E_{t-1} \left [ {\bf n} \right ] +  \sum_{R_j \in R} r_j \Delta n_j.
\end{eqnarray*}
with a solution
\[
E_t \left [ {\bf n} \right ] = (t-k) \sum_{R_j \in \mathcal{R}} r_j \Delta n_j + E_k \left [ {\bf n} \right ]
\]
from which the result follows.

\end{proof}

The expected number of edges can be calculated similarly, one just has to change the elements that one counts. Therefore
\[
E_t \left [ {\bf m} \right ] =  \sum_{i=0}^{\infty}{i p_t \left ({\bf m} = i \right )},
\]
and the expected change in the number of edges can be calculated as
\[
\lim_{t \rightarrow \infty}\frac{E_t \left [ {\bf m} \right ]}{t} = \sum_{R_i \in \mathcal{R}}{r_i \Delta m_i}.
\]

The calculation of expected number of vertices and edges in a step $t$ for a preferential attachment model (see section \ref{s:PA model}) is straightforward. Because of the application of a single rule, we know the exact numbers (the change of vertices and edges in a single time step is 1), and since after the first step $t=1$, we have already 3 vertices and 2 edges, we obtain
\[
E_t \left [ {\bf n}_{pa} \right ] = t+2 \;\;\; \mbox{and} \;\;\;
E_t \left [ {\bf m}_{pa} \right ] = t+1.
\]

% --- intro to rule selection PICGs ---%

\begin{figure*}[thbp]
\setlength{\unitlength}{0.6pt}
  \begin{picture}(500,170)
   \put(0,20){\includegraphics[width=45mm]{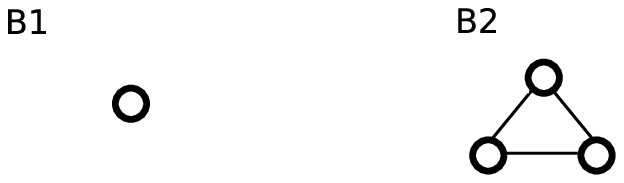}}
   \put(0,0){\small (a) basic graphs}
   \put(250,20){\includegraphics[width=55mm]{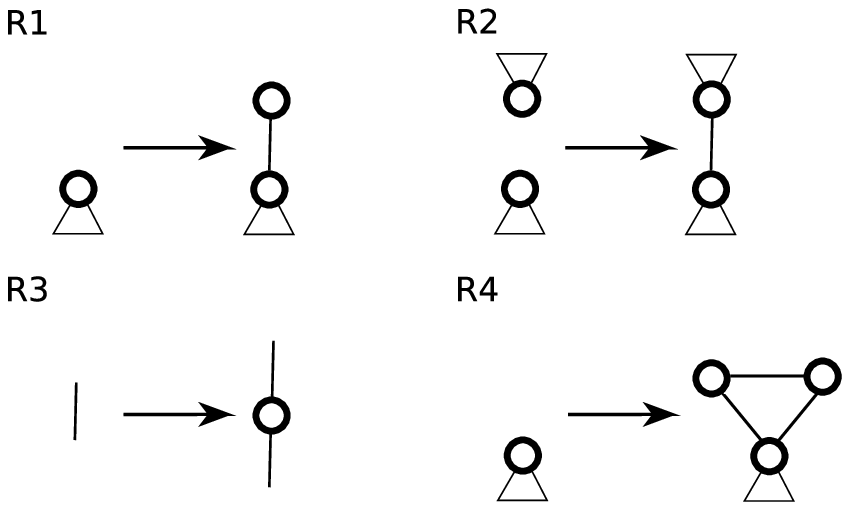}}
   \put(250,0){\small (b) rules}
  \end{picture}
\caption{\label{f:bases_rules}Bases (\textit{B1}, \textit{B2}) and rules for the PICG models presented.}
\end{figure*}

In the next three sections, we present three PICG models. The ICG models are defined by (1) subset of the basis graphs, and (2) the subset of rules from Fig. \ref{f:bases_rules}. Rules are \textit{R1} --- add a vertex with an edge to an arbitrary vertex with degree $i \geq 0$; \textit{R2} --- add an edge to two vertices with arbitrary degrees (the vertices of the left side are different, but can be endpoints of the same edge, thus the rule allows for multiple edges but no loops); \textit{R3} --- expand an edge with a vertex and another edge; \textit{R4} --- add two interconnected vertices to a vertex with arbitrary degree. They share some common characteristics. (3) The basis of all models considered in this paper contains only one graph, therefore it is chosen with probability 1; (4) each of the rules in a model is chosen with a non-zero probability and; (5) the left elements from the rule are chosen uniformly at random.

We then show how one can derive some probabilistic properties of graphs that belongs to the defined classes.

% ----------------------------------------------------------------------------
% --- intro to rule selection PICGs ---%

\section{\label{s:connected}Connected undirected graphs}

The ICG description is $\mathcal{I}(B1;R1,R2)$ according to the rules and the basis from the Fig. \ref{f:bases_rules}. For PICG extension, $\mathcal{I}_{C}$, we add the probability for rule selection: \textit{R1} is selected with probability $q$ and \textit{R2} with probability  $r=1-q$, $q \in (0,1)$.

\subsection{Order and size}
\label{ss:connected order_size}

In the very precise formulation of the model we note, that in the first step only \textit{R1} can be applied, but the initial condition does not depend on the further evolution of the graph.

The number of edges at step $t$ is $m$ if and only if the number of edges at step $\left ( t-1 \right )$ is $\left ( m - 1 \right )$ and any rule is applied, therefore:

\begin{eqnarray*}
p_t \left ( {\bf m} = m \right ) &=& q p_{t-1} \left ( {\bf m} = m - 1\right ) + r p_{t-1} \left ( {\bf m} = m-1 \right ) =\\
&=& p_{t-1} \left ( {\bf m} = m-1 \right ).
\end{eqnarray*}

\begin{lem}
The probability that a graph $G \in \mathcal{I}_{C}$ has $m$ edges at step $t$ is given by:

\[
p_t\left ( {\bf m} = m \right ) = \left\{
  \begin{array}{ll}
    1 & t = m \\
    0 & t \ne m\\
  \end{array}. \right.
\]

\end{lem}

\begin{proof}
The lemma follows by applying the initial condition, $p_0 \left ( {\bf m} = 0 \right ) = 1$, to the previously stated recurrence.
\end{proof}

The number of vertices at step $t$ is $n$ if and only if the number of vertices at step $\left ( t - 1 \right )$ is: (1) $n$ and rule \textit{R1} is applied, (2) $\left (n-1 \right )$ and rule \textit{R2} is applied, or (3) $n$ and neither rule $R1$ nor rule $R2$ is applied (thus, the probability is 0).

\begin{lem}
The probability that a graph $G \in \mathcal{I}_{C}$ has $n$ nodes at step $t$ is given by:

\[
p_t\left ( n \right ) = \left\{
  \begin{array}{ll}
    {{t-1} \choose {t - n + 1}} q^{n-2} r^{t-n+1} & n > 1 \\
    0 & n = 0 \; or \; 1\\
  \end{array}. \right.
\]

\end{lem}

\begin{proof}
For the proof we have that:

\[
p_t \left ( n \right ) = q p_{t-1} \left ( n-1 \right ) + r p_{t-1} \left ( n \right ).
\]

Let $P_n \left (x \right ) = \sum_{t}{p_t \left ( n \right )x^t}$. Multiplying the last expression by $x^t$ and summing over all $t$, we get:

\[
P_n \left ( x \right ) = xrP_n \left ( x \right ) + x q P_{n-1} \left ( x \right ),
\] or equivalently,

\[
P_n \left ( x \right ) = \frac{xq}{1 - x r}P_{n-1} \left ( x \right ).
\]

To determine the initial condition, \textit{i.e.} $P_{2} \left ( x \right )$, we find that:

\[
p_t \left ( 2 \right ) = rp_{t-1} \left ( 2 \right ) = r^{t-1},
\] since $p_1 \left ( 2 \right ) = 1$ and for every $t > 0$, $p_t \left ( {\bf n} = n \right ) = 0$ when $n = 0$ or $1$. Therefore, $P_{2} \left ( x \right ) = \frac{x}{1- rx}$. To obtain the desired result, we observe that $p_t \left ( n \right )$ is the coefficient of the $x^t$ term in the expansion of $P_n \left (x \right )$:

\begin{eqnarray*}
p_t \left ( n \right ) & = & \left [ x^t \right ] P_n \left ( x \right ) = \\
&=& \left [ x^t \right ] \left ( \frac{xq}{1 - xr}\right )^{n-2}\frac{x}{1-xr} = \\
&=& q^{n-2}\left [x^{t-n+1} \right ] \frac{1}{\left ( 1 - xr \right )^{\left (n-2 \right) + 1}} = \\
&=& q^{n-2} {{\left (t-n+1\right ) + \left ( n - 2 \right )} \choose {t - n + 1}}  r^{t-n+1} =\\
&=& {{t-1} \choose {t - n + 1}} q^{n-2} r^{t-n+1},
\end{eqnarray*}
and we have the claim.

\end{proof}

Finally, we determine the expected number of nodes at step $t$. The result is summarized in the following theorem and corollary:

\begin{thm}
Given a graph $G \in \mathcal{I}_{C}$, the expected number of nodes at step $t$, denoted by $E_t \left [ {\bf n} \right ]$, is:

\[
E_t \left [ {\bf n} \right ] = \left ( t - 1 \right )q + 2 - \left ( t + 1 \right ) q^{t-1}.
\]
\end{thm}

\begin{proof}
The proof involves some technical manipulation of the expression for the expectation:

\begin{eqnarray*}
E_t \left [ {\bf n} \right ] &=& \sum_{i = 2}^{t}{i {{t-1} \choose {t - i + 2}} q^{i-2} r^{t-i+1}} =\\
&=& \sum_{j = 0}^{t-2}{\left (j+2 \right ) {{t-1} \choose {t - j - 1}} q^{j} r^{t-j-1}} =\\
&=& \sum_{j = 0}^{t-1}{\left (j+2 \right ) {{t-1} \choose {t - j - 1}} q^{j} r^{t-j-1}} - \left ( t + 1 \right ) q^{t-1}=\\
&=& \sum_{j = 0}^{t-1}{j {{t-1} \choose {t - j - 1}} q^{j} r^{t-j-1}} + 2 - \left ( t + 1 \right ) q^{t-1}.
\end{eqnarray*}

Since $\sum_{j = 0}^{t}{j {{t} \choose {t - j}} q^{j}r^{t-j}} = tq$ (the expected value of the binomial distribution), we have the proof.

\end{proof}

\begin{cor}
When $t \rightarrow \infty$, $\frac{E_t \left [ {\bf n} \right ]}{t} \rightarrow q$ in probability.
\end{cor}

\begin{proof}
Finding the limit $\frac{E_t \left [ {\bf n} \right ]}{t}$ when $t \rightarrow \infty$ is an elementary exercise.
\end{proof}

\subsection{\label{ss:deg_distr_1}Degree distribution}

\begin{figure}[!thbp]
\begin{center}
 \includegraphics[width=80mm]{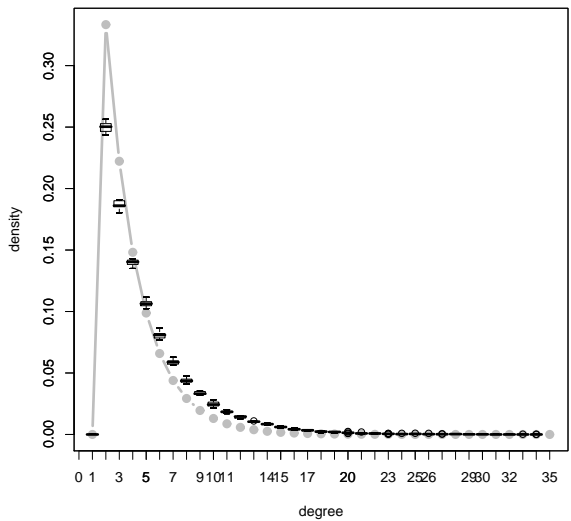}
\end{center}
\caption{\label{f:degree fittings C}10 realizations of the model were grown for $q = r = 0.5$ till they reached $10,000$ vertices. The degree density was calculated. Boxplot for each degree is represented on the graph (in black). We can see that variance of the degree densities is very small. Gray points represent the the mean-field degree distribution fit. The lines are there just for easier visual inspection.}
\end{figure}

Degree distribution of several growing graphs can be approximately calculated with the use of generating functions \cite{Moore06}. Consider $p_d$ the probability of a vertex having a degree $d$. If there are $n$ vertices currently in the graph, there are on average $np_d$ vertices of degree $d$ currently in the graph.
With this we assume, that the number of vertices of degree $d$ is tightly concentrated around its expectation ($np_d$), which corresponds with the simulation results for large $n$ (see Fig. \ref{f:degree fittings C}). A claim similar to this one can be found rigourosly proven in the paper of Cooper et al \cite{Cooper04} for  scale free random graph processes.

In the next step, the probability $p_d$ is going to be updated regarding the use of a rule. We denote it as $p_d'$. The number of all vertices in the next step will be larger for the probability, that a new vertex is added, so $n + q$ (only \textit{R1} adds a vertex). And the number of vertices with degree $d$ present in the next step will be $(n + q)p_d'$.

In this way, we can write the rate equation for the evolution of the graph degree distribution. The number of vertices with degree $d$ changes if either of the rules is applied. With probability $q$, the \textit{R1} is used, a vertex with degree 1 is added (linked to a randomly selected vertex). In the equation we describe the addition of this vertex by Kronecker $\delta_{d1}$. The choice of a vertex can influence $p_d$ iff the chosen vertex is of degree $d$ ($p_d$ is diminished) or of degree $(d-1)$ ($p_d$ increases).

\textit{R2} is applied with probability $r$.
The vertices for rule \textit{R2} are chosen uniformly at random. When a vertex with degree $(d-1)$ is chosen, this adds to the number of vertices with degree $d$, but when a vertex with degree $d$ is chosen, the number of vertices with degree $d$ will be diminished. The probability, that two vertices of degree $d$ are selected is
\[
\frac{\frac{np_d(np_d-1)}{2}}{\frac{n(n-1)}{2}} = \frac{p_d(np_d-1)}{n-1} \approx p_d^2.
\label{eq:p_d contrib}
\]
The entire change in the number of vertices with degree $d$ in one time step can therefore be written as
\[
  2p_{d-1}^2 - 2p_d^2 + p_{d-1}(1-p_{d-1}-p_d) + p_{d}(1-p_{d-1}-p_d) = p_{d-1} - p_{d} + p_{d-1}^2 - p_d^2.
\]
The contributions, where both vertices change the $p_d$ probability simultaneously vanish in the limit of large $n$ and have therefore been neglected.
The rate equation follows:
\[
(n+q)p'_d = np_d + q\left (\delta_{d1} + p_{d-1} - p_d \right ) + r\left (p_{d-1} - p_d\right ).
\]
$p'_d$ is taken to equal $p_d$, when $n \rightarrow \infty$, which gives
\[
(1+q)p_d - q\delta_{d1} - p_{d-1} = 0.
\]
The solution of this equation can be found in terms of generating functions. Multiplying by $z^d$ and summing over $d$ (considering $\sum_d p_d z^d = G(z)$)
\begin{eqnarray*}
(1+q)G(z) - qz - z G(z) &=& 0 \\
G(z) &=& \frac{qz}{(1+q)-z}.
\end{eqnarray*}
The $[z^i]$-th term of the generating function $G(z)$ corresponds to the probability $p_i$ and rewriting the generating function in terms of geometric series we get
\[
 [z^i]g(z) = \frac{G^{(i)}(z)}{i!}\big |_{z=0} = \frac{q}{(1+q)^{i}}.
\]
This shows, that the degree distribution of the graph follows the exponential function.

Considering the rule \textit{R2} only, where vertices in the left element are not linked, the inductive description $\mathcal{I}(B1;R1,R2)$ represents exactly the class of all \textit{simple connected undirected graphs}. Note, that in this case, the left elements for \textit{R2} are not selected uniformly at random anymore. However this anomaly disappears, when $t \rightarrow \infty$, while the probability, that the two random selected vertices share an edge goes to 0.

The model above is therefore a good approximation for the class of all simple connected undirected graphs. This can also be seen from the simulation of the degree distribution (see Fig. \ref{f:degree fittings C}).

% --------------------------------------- %

\section{2-vertex-connected graphs}

A graph $G$ is \textit{2-connected} if $|G| > 2$ and $G - u$ is connected for every $u \in V(G)$. In other
words, no pair of vertices can be separated by the removal of any other vertex \cite{Diestel05}.
The ICG description is $\mathcal{I}(B2;R2,R3)$ according to the rules and the bases defined above \cite{Lovasz79}. In the probabilistic extension $\mathcal{I}_{2V}$, \textit{R2} is selected with probability $q \in (0,1)$, and \textit{R3} with probability $r$, where $q+r=1$.

\subsection{Order and size}

The calculations for order and size are done very similarly to the ones in subsection \ref{ss:connected order_size}. The change in the number of vertices and in the number of edges with the application of the rules is the same. The number of edges changes for 1 in every step and the number of vertices does not change when \textit{R2} is applied (with probability $q$) and it changes for one when \textit{R3} is applied (with probability $r$). Note that the probabilities of choosing a rule are opposite to the ones in the subsection \ref{ss:connected order_size}. The second difference comes with the different initial graph; \textit{B2} has already 3 vertices and 3 edges.

\begin{lem}
The probability that a graph $G \in \mathcal{I}_{2V}$ has $m$ edges at step $t$ is given by:

\[
p_t\left ( {\bf m} = m \right ) = \left\{
  \begin{array}{ll}
    1 & t = m+3 \\
    0 & t \ne m+3\\
  \end{array}. \right.
\]

\end{lem}

\begin{proof}
The lemma follows by applying the initial condition, $p_0 \left ( {\bf m} = 3 \right ) = 1$, to the recurrence stated in subsection \ref{ss:connected order_size}.
\end{proof}

\begin{lem}
The probability that a graph $G \in \mathcal{I}_{2V}$ has $n$ nodes at step $t$ is given by:

\[
p_t\left ( n \right ) = \left\{
  \begin{array}{ll}
    {{t} \choose {t - n + 3}} r^{n-3} q^{t-n+3} & n > 2 \\
    0 & n \in \{0,1,2\}\\
  \end{array}. \right.
\]

\end{lem}

\begin{proof}
Similarly to the proof in subsection \ref{ss:connected order_size}:

\[
p_t \left ( n \right ) =q p_{t-1} \left ( n \right ) + r p_{t-1} \left ( n-1 \right ).
\]

For the generating function's expression, when $P_n \left (x \right ) = \sum_{t}{p_t \left ( n \right )x^t}$ we get:

\[
P_n \left ( x \right ) = xqP_n \left ( x \right ) + x r P_{n-1} \left ( x \right )
= \frac{xr}{1 - x q}P_{n-1} \left ( x \right ).
\]

To determine the initial condition, \textit{i.e.} $P_{3} \left ( x \right )$, we find that $ p_t \left ( 3 \right ) = qp_{t-1} \left ( 3 \right ) = q^{t}$, since $p_0 \left ( 3 \right ) = 1$ and for every $t > 0$, $p_t \left ( n \right ) = 0$ when $n \in \{0,1,2\}$. Therefore, $P_{3} \left ( x \right ) = \frac{1}{1- qx}$.

Again, $p_t \left ( n \right )$ is the coefficient of the $x^t$ term in the expansion of $P_n \left (x \right )$:
\[
p_t \left ( n \right ) = {{t} \choose {t - n + 3}} r^{n-3} q^{t-n+3},
\]
and we have the claim.

\end{proof}

The expected number of nodes at step $t$ is summarized in the following theorem:

\begin{thm}
Given a graph $G \in \mathcal{I}_{2V}$, the expected number of nodes at step $t$, denoted by $E_t \left [ {\bf n} \right ]$, is:

\[
E_t \left [ {\bf n} \right ] = tr + 3 - \frac{(t+1)t(t-1)}{2}q^2r^{t-2} - (t+2)t qr^{t-1} - (t-3)r^t.
\]
\end{thm}

\begin{proof}
The proof involves some technical manipulation of the expression for the expectation:

\begin{eqnarray*}
E_t \left [ {\bf n} \right ] &=& \sum_{i = 3}^{t}{i {{t} \choose {t - i + 3}} q^{t-i+3} r^{i-3}} = \sum_{j = 0}^{t-3}{\left (j+3 \right ) {{t} \choose {t - j}} q^{t-j} r^{j}} =\\
&=& \sum_{j = 0}^{t}{\left (j+3 \right ) {{t} \choose {t - j}} q^{t-j} r^{j}} - \frac{(t+1)t(t-1)}{2}q^2r^{t-2} - (t+2)t qr^{t-1} - (t-3)r^t=\\
&=& \sum_{j = 0}^{t}{j {{t} \choose {t - j}} q^{t-j} r^{j}} + 3 - \frac{(t+1)t(t-1)}{2}q^2r^{t-2} - (t+2)t qr^{t-1} - (t-3)r^t.
\end{eqnarray*}

Since $\sum_{j = 0}^{t}{j {{t} \choose {t - j}} r^{j}q^{t-j}} = tr$ (the expected value of the binomial distribution), we have the proof.

\end{proof}

\begin{cor}
When $t \rightarrow \infty$, $\frac{E_t \left [ {\bf n} \right ]}{t} \rightarrow r$ in probability.
\end{cor}

\subsection{Degree distribution.}

We use the approximate, mean-field approach (as described in \ref{ss:deg_distr_1}), solved by the generating functions methodology. The rate equation for the change in the number of vertices of degree $d$ is written in the following way
\[
 (n+r)p'_d = np_d + r\left (\delta_{d2} \right ) + q\left (p_{d-1} - p_d\right ).
\]
$p'_d$ is taken to equal $p_d$, when $n \rightarrow \infty$, which gives
\[
 p_d - r\delta_{d2} - qp_{d-1} = 0.
\]
Multiplying by $z^{d-2}$ and summing over $d$ (considering $\sum_{d \geq 2} p_d z^{d-2} = F(z)$ --- left shifted generating function, where $p_0 = 0$ and $p_1 = 0$)
\begin{eqnarray*}
 F(z) - r - qz F(z) &=& 0 \\
 F(z) &=& \frac{r}{1-qz}.
\end{eqnarray*}
The $[z^{i-2}]$-th term ($i \geq 2$) of the generating function $F(z)$ corresponds to the probability $p_i$,
\[
 [z^{i-2}]F(z) = \frac{F^{(i-2)}(z)}{(i-2)!}\big |_{z=0} = rq^{i-2}.
\]

\begin{figure}[!thbp]
 \includegraphics[width=80mm]{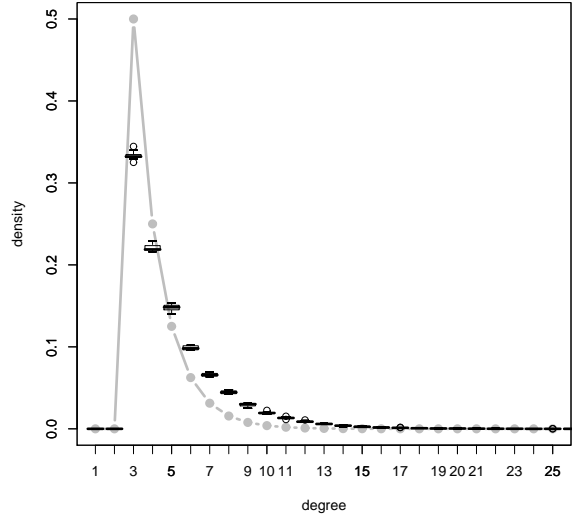}
\caption{\label{f:degree fittings 2V}10 realizations of the model were grown for $q = r = 0.5$ till they reached $10,000$ vertices. Boxplot for the density of each degree is represented on the graph (in black). Gray points represent the the mean-field degree distribution fit. The lines are there for easier visual inspection.}
\end{figure}

Figure \ref{f:degree fittings 2V} shows the fittings of degree distribution to the simulated data.

% --------------------------------------- %

\section{2-edge-connected graphs}

$G$ is called \textit{2-edge-connected} if $|G| > 1$ and $G - e$ is connected for every $e \in E(G)$ \cite{Diestel05}.
The ICG description is $\mathcal{I}(B2;R2,R3,R4)$ according to the rules and the basis defined above \cite{Lovasz79}. For the PICG model $\mathcal{I}_{2E}$, \textit{R2} is selected with probability $q$, \textit{R3} with probability $r$ and \textit{R4} with probability $s$, $q+r+s=1$.

%%%% --------------------------------------------------------------- %
\subsection{Order and size}

The number of vertices at step $t$ is $n$, $n > 4$, if and only if the number of vertices at step $(t-1)$ is: (1) $n$ and rule \textit{R2} is applied, or (2) $n-1$ and rule \textit{R3} is applied, or (3) $(n-2)$ and rule \textit{R4} is applied:
\[
p_t(n) = q p_{t-1}(n) + r p_{t-1}(n-1) + s p_{t-1}(n-2).
\]

Consider $P_n(x)=\sum_{t\ge 0} {p_t(n)x^t}$. Multiplying the last expression and summing over $t \ge 0$ gives

\[
P_n(x) = q x P_n(x) + r x P_{n-1}(x) + s x P_{n-2}(x),
\] or equivalently:

\[
(1-qx)P_n(x) = rxP_{n-1}(x) + sxP_{n-2}(x).
\]

Note that $P_0(x) = P_1(x) = P_2(x) = 0$. From the initial condition, we further have that $p_0(3) = 1$, and $p_t(3) = q^tp_0(3) = q^t$, therefore, $P_3(x) = \frac{1}{1-qx}$.

The following claim can be inductively proven:

\begin{lem}
Given a positive integer $n$, the following holds for the generating function $P_n(x)$:
\begin{itemize}
\item $P_0(x) = P_1(x) = P_2(x) = 0$,
\item $P_3(x) = \frac{1}{1-qx}$,
\item For every $n \ge 3$, $P_n(x)$ is a summation of $\lceil \frac{n}{2} \rceil - 1$ functions:

\[
P_n(x) = \sum_{j = \lfloor \frac{n}{2} \rfloor}^{n-2}{{{2j-n+2} \choose {n-j-2}} \frac{r^{2j - n+1}s^{n-j-2}x^{j-1}}{(1 - qx)^j}}.
\]
\end{itemize}
\end{lem}

Then we have the following theorem:

\begin{thm}
The probability that a 2-edge connected random graph has $n$ vertices after applying $t$ rules is given by:

\[
p_t(n) = \sum_{j = \lfloor \frac{n}{2} \rfloor}^{n-2}{{{2j-n+2} \choose {n-j-2}} {{t} \choose {t-j+1}}q^{t-j+1}r^{2j - n+1}s^{n-j-2}}.
\]
\end{thm}

\begin{proof}
Let $a_j^n = {{2j-n+2} \choose {n-j-2}}$.
Note that

\begin{eqnarray*}
p_t(n) &=& [x^t] P_n(x) = \\
&=& [x^t]\sum_{j = \lfloor \frac{n}{2} \rfloor}^{n-2}{a_j^n\frac{r^{2j - n+1}s^{n-j-2}x^{j-1}}{(1 - qx)^j}}=\\
&=&\sum_{j = \lfloor \frac{n}{2} \rfloor}^{n-2}{a_j^nr^{2j - n+1}s^{n-j-2}[x^{t-j+1}]\frac{1}{(1 - qx)^j}} = \\
&=& \sum_{j =\lfloor \frac{n}{2} \rfloor}^{n-2}{a_j^nr^{2j - n+1}s^{n-j-2}{{t} \choose {t-j+1}}q^{t-j+1}} = \\
&=& \sum_{j = \lfloor \frac{n}{2} \rfloor}^{n-2}{a_j^n{{t} \choose {t-j+1}}q^{t-j+1}r^{2j - n+1}s^{n-j-2}}.
\end{eqnarray*}

\end{proof}

According to the approximate calculation of the expected change in the number of vertices (Eq. (\ref{eq:expected_change})), the expected change for 2-edge-connected graph model is $
\frac{E_t \left [ {\bf n} \right ]}{t} \approx r+2s$.

%%%% ----------------------------------------------------------------- %

The calculation of the number of edges involves even more demanding and tedious calculations. For the clarity of the presented paper, we do not include that. According to Eq. (\ref{eq:expected_change}) the expected change in the number of edges, when $t \rightarrow \infty$, approximately equals to $(q+r)+3s$.

%%%% ----------------------------------------------------------------- %

\subsection{Degree distribution.}

The mean-field approach gives the rate equation for the change in the number of vertices of degree $d$
\[
 (n+r+2s)p'_d = np_d + r\left (\delta_{d2} \right ) + q\left (p_{d-1} - p_d\right ) + s\left (p_{d-2} + 2\delta_{d2} - p_d \right ).
\]
$p'_d$ is taken to equal $p_d$, when $n \rightarrow \infty$, which gives
\[
 (1+2s)p_d - (r+2s)\delta_{d2} - qp_{d-1} - sp_{d-2} = 0.
\]
Multiplying by $z^{d-2}$ and summing over $d$ (considering $\sum_{d \geq 2} p_d z^{d-2} = F(z)$, again $p_0 = 0$ and $p_1 = 0$)
\begin{eqnarray*}
 & (1+2s)F(z) - (r+2s) - qzF(z) - sz^2F(z) = 0 \\
 & F(z) = \frac{r+2s}{(1+2s)-qz-sz^2}.
\end{eqnarray*}
The $[z^{i-2}]$-th term of the generating function $F(z)$ corresponds to the probability $p_i$. It can be expressed in a closed form according to Graham et al. \cite{Graham} with the partial fraction expansion.

We can write $F(z) = P(z)/Q(z)$ and express $Q(z) = (1+2s)(1-\rho_1z)(1-\rho_2z)$, because the coefficients of the generating functions of the form $1/(1-az)$ are terms of simple geometric series expansion. Calculate the two roots of the ``reflected" polynomial $Q^R(z)$ ($\rho_{1,2}$), that represent the roots of the above expression of $Q(z)$. $[z^{n}]$-th term of $F(z)$ can then be expressed as $a_1 \rho_1^n + a_2 \rho_2^n$ (when the roots are distinct). $a_i$ is expressed as $-\rho_i(P(1/\rho_i))/(Q'(1/\rho_i))$.

In our case, the roots of $Q^R(z)$ equal
\[
 \rho_{1,2} = \frac{q \pm \sqrt{q^2 + 4s(1+2s)}}{2(1+2s)}
\]

The $[z^{i}]$-th coefficient of $F(z)$ therefore is

\[
 [z^i]F(z) = a_1 \cdot \left ( \frac{q + \sqrt{q^2+4s(1+2s)}}{2(1+2s)}\right )^i + a_2 \cdot \left ( \frac{q - \sqrt{q^2+4s(1+2s)}}{2(1+2s)} \right )^i,
\]
where
\begin{eqnarray*}
 a_1 &=& \frac{(q + \sqrt{q^2 + 4s(1+2s)})(r+2s)}{2(1+2s) \sqrt{q^2 + 4s(1+2s)}} \;\;\;\;\; {\text{and}}\\
 a_2 &=& - \frac{(q - \sqrt{q^2 + 4s(1+2s)})(r+2s)}{2(1+2s) \sqrt{q^2 + 4s(1+2s)}}
\end{eqnarray*}

\begin{figure}[!thbp]
 \includegraphics[width=80mm]{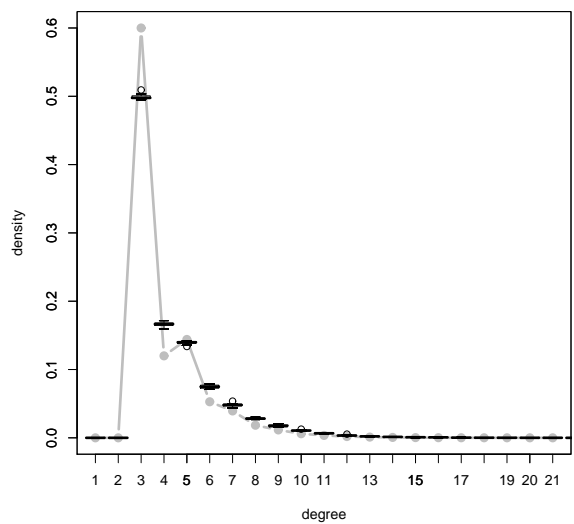}
\caption{\label{f:degree fittings 2E}10 realizations of the model were grown for $q = r = s = 1/3$ till they reached $10,000$ vertices. Boxplot for the density of each degree is represented on the graph (in black). Gray points represent the the mean-field degree distribution fit. The lines are there for easier visual inspection.}
\end{figure}

Figure \ref{f:degree fittings 2E} shows the fittings of degree distribution to the simulated data.

\section{Conclusion}

In the paper we proposed a quite broad framework for the description of statistically evolving complex networks, the probabilistic inductive classes of graphs. We presented, how the already existing network models can be described in this way. We give a general recursive equation for order and size of the classes of graphs, which belong to the simple PICGs and apply them on three PICG models of growing networks.

% ------------------------------------------------------------------------
%Included for Gather Purpose only:
%input "XbibNatasa.bib"
\bibliographystyle{amsplain}
\bibliography{PICGs.bib}
\end{document}